\documentclass[a4paper,12pt]{article}
\usepackage{amsfonts}
\usepackage[xdvi]{graphicx}
\usepackage[dvips]{color}
\usepackage{amsmath,amssymb,amsthm}

\newcommand{\R}{{\mathbb R}}
\newcommand{\C}{{\mathbb C}}

\makeatletter

\@addtoreset{equation}{section}

\newcounter{def}[section]
\renewcommand{\thedef}{\stepcounter{def}\thesection.\@arabic\c@def }

\begin{document}
\setlength{\baselineskip}{24pt}
\begin{center}
\textbf{\LARGE{A center manifold reduction of the Kuramoto-Daido model with a phase-lag}}
\end{center}

\setlength{\baselineskip}{14pt}

\begin{center}
Institute of Mathematics for Industry, Kyushu University / JST PRESTO,\\
Fukuoka, 819-0395, Japan 
\end{center}
\begin{center}
Hayato CHIBA\footnote{chiba@imi.kyushu-u.ac.jp}
\end{center}
 \begin{center}
Sep 11, 2016; Revised Dec 15, 2016
 \end{center}

\begin{center}
\textbf{Abstract}
\end{center}

A bifurcation from the incoherent state to the partially synchronized state of the 
Kuramoto-Daido model with the coupling function $f(\theta ) = \sin (\theta +\alpha _1)
 + h\sin 2(\theta +\alpha _2)$ is investigated based on the generalized spectral theory and 
the center manifold reduction.
The dynamical system of the order parameter on a center manifold is derived
under the assumption that there exists a center manifold on the dual space
of a certain test function space.
It is shown that the incoherent state loses the stability at a critical coupling strength
$K=K_c$, and a stable rotating partially synchronized state appears for $K>K_c$.
The velocity of the rotating state is different from the average of natural frequencies
of oscillators when $\alpha _1 \neq 0$.


\section{Introduction}

\begin{figure}
\begin{center}
\includegraphics[]{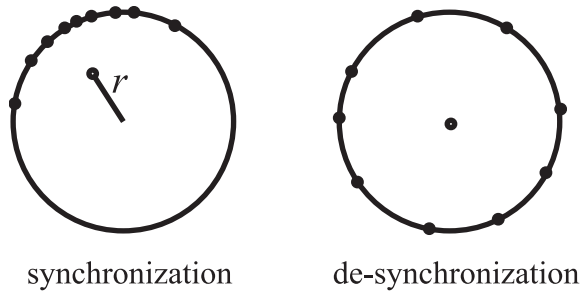}
\caption[]{Collective behavior of oscillators.}
\label{fig1}
\end{center}
\end{figure}

Collective synchronization phenomena are observed in a variety of areas such as chemical reactions,
engineering circuits and biological populations~\cite{Pik}.
In order to investigate such phenomena, a system of globally coupled phase oscillators called the
Kuramoto-Daido model \cite{Dai}
\begin{equation}
\frac{d\theta _i}{dt} 
= \omega _i + \frac{K}{N} \sum^N_{j=1} f (\theta _j - \theta _i),\,\, i= 1, \cdots  ,N,
\label{KM}
\end{equation}
is often used,
where $\theta _i = \theta _i(t) \in [ 0, 2\pi )$ is a dependent variable 
which denotes the phase of an $i$-th oscillator on a circle,
$\omega _i\in \R$ denotes its natural frequency drawn from some density function $g(\omega )$, $K>0$ is a coupling strength,
and where $f(\theta )$ is a $2\pi$-periodic function.
The order parameter defined by
\begin{eqnarray}
r := \Bigl| \frac{1}{N}\sum^N_{j=1}e^{i\theta _j (t)} \Bigr|, \quad i = \sqrt{-1}
\label{1-2}
\end{eqnarray}
is used to measure the amount of coherence in the system;
if $r$ is nearly equal to zero, oscillators are uniformly distributed (called the incoherent state), while if $r>0$,
the synchronization occurs, see Fig. \ref{fig1}.
In the last few decades, the existence, stability and bifurcations of the synchronization
have been well studied, however, the coupling function $f(\theta )$ was still restricted to some specific form.

In this paper, the onset of the synchronization of the continuous limit (thermodynamics limit) for the following model
\begin{equation}
\frac{d\theta _i}{dt} 
= \omega _i + \frac{K}{N} \sum^N_{j=1}\bigl( \sin (\theta _j - \theta _i + \alpha _1)
 + h\cdot \sin 2(\theta_j - \theta _i  + \alpha _2) \bigr) ,
\label{KMP}
\end{equation}
will be considered, where $\alpha _1, \alpha _2$ are phase lags and $h$ is a parameter which 
controls the strength of the second harmonic.
For the continuous limit of the system given in Sec.2, a bifurcation from the incoherent state
(the de-synchronous state) to the partially synchronized state will be investigated based on
the generalized spectral theory.

When $h = \alpha _1 = 0$, the model is well known as the original Kuramoto model.
For a bifurcation of the partially synchronized steady state in this case, it has been proved (the Kuramoto conjecture) that :
\\[0.2cm]
Suppose that the density function $g(\omega )$ for natural frequencies is an even and unimodal function,
and satisfies a certain regularity condition.
\\
\textbf{(i) (Instability of the incoherent state).}
When $K > K_{c}:=2/(\pi g(0))$, then the incoherent state of the continuous limit is unstable.
\\
\textbf{(ii) (Local stability of the incoherent state).}
When $0<K<K_{c}$, the incoherent state is locally asymptotically stable with respect to a certain topology.
\\
\textbf{(iii) (Bifurcation).}
There exists a positive constant $\varepsilon_0$ such that
if $K_{c}-\varepsilon _0 < K < K_{c} + \varepsilon_0 $ and if an initial condition is close to the incoherent state
(with respect to a certain topology), 
the order parameter is locally governed by the dynamical system
\begin{eqnarray}
\frac{dr}{dt} = \text{const.} \left(K-K_c + \frac{\pi g''(0) K_c^4}{16} r^2 \right) r+ O(r^4).
\label{1-4}
\end{eqnarray}
In particular, when $K>K_c$, the order parameter tends to the positive constant expressed as
\begin{equation}
r= \sqrt{\frac{-16}{\pi K_{c}^4 g''(0)}}\sqrt{K - K_{c}} + O(K-K_{c}), 
\end{equation}
as $t\to \infty$.

In Chiba \cite{Chi1, Chi2}, this result is proved based on the generalized spectral theory \cite{Chi3}
under the assumption that $g(\omega )$ has an analytic continuation near the real axis.
In Dietert \cite{Die} and Fernandez et al.~\cite{Fer}, they proved a similar result with a weaker 
assumption for $g(\omega )$.
See these references for the precise assumptions for $g(\omega )$ and for the topologies for the stability.

\begin{figure}
\begin{center}
\includegraphics[scale=1.3]{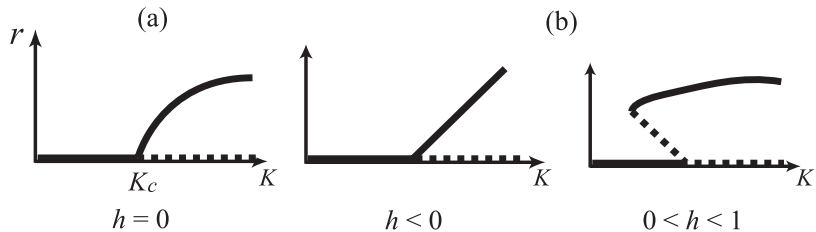}
\caption{Bifurcation diagrams of the order parameter for (a) $f(\theta ) = \sin \theta $
and (b) $f(\theta ) = \sin \theta  + h\sin 2\theta$.
The solid lines denote stable solutions, and the dotted lines denote unstable solutions. \label{fig2}}
\end{center}
\end{figure}

In \cite{Chi1}, this result was extended to the case $\alpha _1 = \alpha _2 = 0$ and $h < 1$ 
(i.e. the coupling function is $f(\theta ) = \sin \theta + h\sin 2\theta $).
With the aid of the generalized spectral theory, the dynamics of the order parameter on a center manifold
was shown to be 
\begin{eqnarray*}
\frac{dr}{dt} = \text{const.} \left(K-K_c - \frac{K_c^2 Ch}{1-h}r \right) r + O(r^3),
\end{eqnarray*}
where $C$ is a certain negative constant.
As a result, a bifurcation diagram of $r$ is given as Fig.\ref{fig2}.
When $h=0$, the synchronous state bifurcates through the pitchfork bifurcation,
though when $h\neq 0$, it is a transcritical bifurcation.

Omel'chenko et al.\cite{OM} have investigated a bifurcation diagram for the case $h=0,\, \alpha _1 \neq 0$
by the self-consistent approach and numerical simulation.
They found that under the synchronized state, the average velocity of oscillators is different from
the average of their natural frequencies because of the phase-lag $\alpha _1$.

The purpose in this paper is to derive the dynamics of the order parameter 
for the system (\ref{KMP}) for a certain parameter region $(\alpha _1, \alpha _2, h)$.
To this end, we need four assumptions (A1) to (A4) given after Section 3.
Here, we give a rough explanation of these assumptions.
\\[0.2cm]
\textbf{(A1)} We assume that $\cos \alpha _1 > 0$ so that the interaction between oscillators through 
the coupling $\sin (\theta +\alpha_1)$ is an attractive coupling.
Further, we assume that $\sin (\theta +\alpha _1) $ in the coupling function is stronger than the 
second harmonic $h\sin 2(\theta + \alpha _2)$ in a certain sense. When $\alpha _1 = 2\alpha _2$, this is equivalent to $h<1$.
\\
\textbf{(A2)} We assume that the density $g(\omega )$ of natural frequencies is an analytic function
near the real axis. This is the essential assumption to apply the generalized spectral theory.
\\
\textbf{(A3)} We will show that at a bifurcation value $K=K_c$, one of the generalized eigenvalues $\lambda _c$
of a certain linear operator obtained by the linearization of the system lies on 
the point $iy_c$ on the imaginary axis.
We assume that $\lambda _c$ is a simple eigenvalue.
\\
\textbf{(A4)} We assume that as $K$ increases, $\lambda _c$ transversally gets across the imaginary axis
at the point $iy_c$ from the left to the right.
\\

See after Section 3 for the precise statements for these assumptions.
If $g(\omega )$ is an even and unimodal function, (A3) is automatically satisfied,
though we do not assume it in this paper.
The main results in the present paper are;
\\[0.2cm]
\textbf{Theorem \thedef \,(Instability of the incoherent state).}
\\
Suppose (A1) and $g(\omega )$ is continuous.
There exists a number $\varepsilon >0$ such that
when $K_c<K<K_c+\varepsilon $, the incoherent state is linearly unstable
(the value of $K_c$ will be given in Section 3).
\\[0.2cm]
\textbf{Theorem \thedef\, (Local stability of the incoherent state).}
\\
Suppose (A1) and (A2).
When $0<K<K_c$, the incoherent state is linearly asymptotically stable in the weak sense
(see Section 4 for the weak stability).
\\

For a bifurcation, our result is divided into two cases, $h=0$ and $h\neq 0$ because 
types of bifurcations of them are different.
As is shown in Fig.\ref{fig2}, it is a pitchfork bifurcation when $h=0$, while it is a transcritical bifurcation when $h\neq 0$.
In this paper, we formally apply the center manifold reduction without a proof of the existence of a center manifold.
In Chiba \cite{Chi2}, the existence of a center manifold was shown when $h=\alpha _1=0$ and $g$ is the Gaussian.
To prove the existence of center manifolds for a wide class of evolution equations
within the setting of the generalized spectral theory is an important challenging work.
\\[0.2cm]
\textbf{Theorem \thedef\, (Bifurcation $h=0$).}
\\
Suppose (A1) to (A4) hold.
There exists a positive constant $\varepsilon_0$ such that
if $K_{c}-\varepsilon _0 < K < K_{c} + \varepsilon_0 $ and if an initial condition is close to the incoherent state, 
the order parameter is locally governed by the dynamical system
\begin{equation}
\frac{dr}{dt}= \mathrm{Re}(p_1)r \left( K-K_c + \frac{\mathrm{Re}(p_3)}{\mathrm{Re}(p_1)}r^2 \right) + O(r^4) ,
\end{equation}
where $p_1$ and $p_3$ are certain complex constants given in Section 6.
The equation has the fixed point
\begin{equation}
r_0 = \sqrt{-\frac{\mathrm{Re}(p_1)}{\mathrm{Re}(p_3)}} \cdot \sqrt{K-K_c} + O(K-K_c).
\end{equation}
The assumption (A4) implies $\mathrm{Re}(p_1)>0$.
Thus, if $\mathrm{Re}(p_3) < 0$, the bifurcation is supercritical; 
the fixed point exists for $K>K_c$ and is stable.
If $\mathrm{Re}(p_3) > 0$, the bifurcation is subcritical; 
the fixed point exists for $K<K_c$ and is unstable.
\\[0.2cm]
\textbf{Theorem \thedef\, (Bifurcation $h\neq 0$).}
\\
Suppose (A1) to (A4) hold.
There exists a positive constant $\varepsilon_0$ such that
if $K_{c}-\varepsilon _0 < K < K_{c} + \varepsilon_0 $ and if an initial condition is close to the incoherent state, 
the order parameter is locally governed by the dynamical system
\begin{equation}
\frac{dr}{dt}= \mathrm{Re}(p_1)r \left( K-K_c + \frac{\mathrm{Re}(p_2)}{\mathrm{Re}(p_1)}r \right) + O(r^3) ,
\end{equation}
where $p_1$ and $p_2$ are certain complex constants given in Section 6.
The equation has the fixed point
\begin{equation}
r_0 = -\frac{\mathrm{Re}(p_1)}{\mathrm{Re}(p_2)} \cdot (K-K_c) + O((K-K_c)^2).
\end{equation}
The assumption (A4) implies $\mathrm{Re}(p_1)>0$.
Thus, if $\mathrm{Re}(p_2) < 0$, the bifurcation is supercritical; 
the fixed point exists for $K>K_c$ and is stable.
If $\mathrm{Re}(p_2) > 0$, the bifurcation is subcritical; 
the fixed point exists for $K<K_c$ and is unstable.
\\[0.2cm]
\textbf{Theorem \thedef\, (The average velocity).}

The complex order parameter (see Eq.(\ref{2-1}) and (\ref{conti})) 
under the partially locked state shown in Theorems 1.3 and 1.4 is given by
\begin{eqnarray*}
\eta_1(t) = r_0e^{i\alpha _1} \cdot e^{i(y_c + O(r_0))t}(1+ O(r_0)),
\end{eqnarray*}
where $r_0$ is a number given in Theorems 1.3 and 1.4, respectively, for $h=0$ and $h\neq 0$.
Hence, the average velocity of locked oscillators is approximately given by $y_c$ (see (A3) and (A4) for the definition of $y_c$).
When $g(\omega )$ is an even and unimodal function and if $\alpha _1$ is small,
$y_c$ is given by
\begin{equation}
y_c = -\frac{g(0)}{H[g]'(0)}\alpha _1 + O(\alpha _1^2),
\end{equation}
where $H[g]$ is a Hilbert transform of $g$, see Lemma 3.1 and Example 3.6.
\\

When $\alpha _1 =\alpha _2=0$, Theorems 1.3 and 1.4 recover the results of \cite{Chi1, Chi2},
and when $h=0$, Theorem 1.3 coincides with the result in \cite{OM}.


\section{The continuous model}

For the finite dimensional Kuramoto-Daido model (\ref{KM}), 
the $k$-th order parameter is defined by
\begin{eqnarray}
\hat{\eta}_k(t) := \frac{1}{N}\sum^N_{j=1} e^{i k\theta _j(t)}.
\label{2-1}
\end{eqnarray} 
By using it, Eq.(\ref{KM}) is rewritten as 
\begin{eqnarray*}
\frac{d\theta _j}{dt} = \omega _j + K \sum^\infty_{l=-\infty} f_l \hat{\eta}_l (t) e^{-i l \theta_j },
\quad f(\theta ) := \sum^\infty_{l=-\infty}f_l e^{il\theta }.
\end{eqnarray*}
This implies that the flow of $\theta _j$ is generated by the vector field
\begin{eqnarray*}
\hat{v} = \omega _j + K \sum^\infty_{l=-\infty} f_l \hat{\eta}_l (t) e^{-i l \theta_j}.
\end{eqnarray*}
Thus, the continuous model of Eq.(\ref{KM}) is the equation of continuity of the form
\begin{equation}
\left\{ \begin{array}{ll}
\displaystyle \frac{\partial \rho}{\partial t} + \frac{\partial }{\partial \theta }(\rho v) = 0,
\quad \rho = \rho (t, \theta , \omega ), 
\label{conti}  \\
\displaystyle v := \omega  + K \sum^\infty_{l=-\infty} f_l \eta_l (t) e^{-i l \theta },  \\
\displaystyle \eta_l(t) := \int_{\R} \! \int^{2\pi}_{0} \!  e^{i l \theta } \rho (t,\theta , \omega ) g(\omega ) d\theta d\omega.
\end{array} \right.
\end{equation}
Here, $g(\omega )$ is a given probability density function for natural frequencies, 
and the unknown function $\rho = \rho (t, \theta , \omega )$ is a probability measure on $[0, 2\pi)$ parameterized by $t, \omega \in \R$.
It is known that a solution of Eq.(\ref{KM}) converges to that of Eq.(\ref{conti}) as $N\to \infty$
in some probabilistic sense \cite{CM}.
Our goal is to investigate the dynamics of $\eta_1 (t)$, which is a continuous version of Kuramoto's order parameter (\ref{1-2}).
Roughly speaking, $\rho (t, \theta , \omega )$ denotes a probability that
an oscillator having a natural frequency $\omega $ is placed at a position $\theta $.
The trivial solution $\rho = 1/(2\pi)$ of the system is a uniform distribution on the circle,
which is called the incoherent state (de-synchronous state).
In this case, $\eta_ l = 0$ for all $l = \pm 1, \pm 2, \cdots $.
We will show that a nontrivial stable solution such that $|\eta _1| > 0$ bifurcates from the incoherent state.

Setting the Fourier coefficients
\begin{eqnarray*}
Z_j(t, \omega ) := \int^{2\pi}_{0} \! e^{ij\theta }  \rho (t, \theta , \omega ) d\theta 
\end{eqnarray*}
yields the system of evolution equations of $Z_j$ as
\begin{eqnarray}
\frac{dZ_j}{dt} = ij\omega Z_j +i j K f_j \eta _j + i j K \sum_{l\neq j} f_l\eta_l Z_{j-l}.
\label{Z}
\end{eqnarray}
Indeed, we have
\begin{eqnarray*}
\frac{dZ_j}{dt}
&=& -\int^{2\pi}_{0}\! e^{ij \theta } \frac{\partial }{\partial \theta }(\rho v) d\theta \\
&=& ij \int^{2\pi}_{0}\! e^{ij\theta } \rho (\omega +K\sum^\infty_{l=-\infty} f_l \eta_l e^{-il\theta }) d\theta  \\
&=& ij\omega Z_j + ijK \sum^\infty_{l=-\infty} f_l \eta_l Z_{j-l}. 
\end{eqnarray*}
Since $Z_0 \equiv 1$ because of the normalization $\int^{2\pi}_{0} \! \rho (t, \theta , \omega ) d\theta =1$,
we obtain Eq.(\ref{Z}).
The trivial solution $Z_j \equiv 0 \, (j=\pm 1, \pm 2, \cdots)$ corresponds to the incoherent state.
In what follows, we consider only the equations for $Z_1, Z_2,\cdots $ because $Z_{-j}$ is the complex
conjugate of $Z_j$.


\section{The transition point formula and linear instability}

To investigate the stability of the incoherent state, we consider the linearized system.
Let $L^2 (\R, g(\omega )d\omega )$ be the weighted Lebesgue space with the inner product
\begin{eqnarray*}
(\phi, \psi) = \int_{\R}\! \phi(\omega )\overline{\psi(\omega )} g(\omega )d\omega . 
\end{eqnarray*}
Put $P_0 (\omega ) \equiv 1 \in L^2 (\R, g(\omega )d\omega )$.
We define the one-dimensional integral operator $\mathcal{P}$ on $L^2 (\R, g(\omega )d\omega )$ to be
\begin{equation}
(\mathcal{P}\phi ) (\omega ) = \int_{\R}\! \phi (\omega ) g(\omega )d\omega  = (\phi, P_0) \cdot P_0(\omega ). 
\end{equation}
Then, the order parameters are written by
\begin{equation}
\eta _j(t) 
= \int_{\R} \! Z_j(t, \omega )g(\omega ) d\omega  = \mathcal{P}Z_j. 
\label{eta _j}
\end{equation}
Note that the term $\eta_l Z_{j-l} = (\mathcal{P}Z_l)Z_{j-l}$ in Eq.(\ref{Z}) is a nonlinear term of $Z_{\pm 1}, Z_{\pm 2},\cdots $ when $j\neq l$.
Hence, the linearized system of Eq.(\ref{Z}) around the incoherent state is given by
\begin{equation}
\frac{dZ_j}{dt} = T_jZ_j := (ij\omega + ijKf_j \mathcal{P}) Z_j, \quad j=1,2,\cdots 
\label{linear}
\end{equation}

Let us consider the spectra of linear operators $T_j$.
The multiplication operator $\phi (\omega ) \mapsto \omega \phi (\omega )$ on
$L^2(\R, g(\omega ) d\omega )$ is self-adjoint.
The spectrum of it consists only of the continuous spectrum given by $\sigma _c(\omega ) = \mathrm{supp} (g)$
(the support of $g$).
Therefore, the spectrum of the multiplication by $ij\omega $ lies on the imaginary axis;
$\sigma _c(ij\omega ) = ij\cdot \mathrm{supp} (g)$
(later we will suppose that $g$ is analytic, so that $\sigma _c(ij\omega )$ is the whole imaginary axis).
Since $\mathcal{P}$ is compact, it follows from the perturbation theory of linear operators \cite{Kato}
that the continuous spectrum of $T_j$ is given by $\sigma _c(T_j) = ij\cdot \mathrm{supp} (g)$,
and the residual spectrum of $T_j$ is empty.

When $f_j \neq 0$, eigenvalues $\lambda$ of $T_j$ are given as roots of the equation
\begin{equation}
\int_{\R} \! \frac{1}{\lambda - ij \omega }g(\omega )d\omega
= \frac{1}{ijKf_j}, \quad \lambda \notin \sigma _c(T_j).
\label{eigen-eq}
\end{equation}
Indeed, the equation $(\lambda -T_j)v = 0$ provides
\begin{eqnarray*}
v = ijKf_j (v,P_0) (\lambda -ij\omega )^{-1} P_0.
\end{eqnarray*}
Taking the inner product with $P_0$, we obtain Eq.(\ref{eigen-eq}).
If $\lambda $ is an eigenvalue of $T_j$, the above equality shows that 
\begin{equation}
v_\lambda (\omega ) = \frac{1}{\lambda -ij\omega }
\label{ef}
\end{equation}
is the associated eigenfunction.
This is not in $L^2(\R, g(\omega )d\omega )$ when $\lambda $ is a purely imaginary number.
Thus, there are no eigenvalues on the imaginary axis.
Putting $\lambda =x+iy$ in Eq.(\ref{eigen-eq}) provides
\begin{eqnarray*}
\left\{ \begin{array}{ll}
\displaystyle \int_{\R}\! \frac{x}{x^2 + (y-j\omega )^2}g(\omega )d\omega
   = \frac{-\mathrm{Im}(f_j)}{jK|f_j|^2},  &  \\[0.4cm]
\displaystyle \int_{\R}\! \frac{y-j\omega }{x^2 + (y-j\omega )^2}g(\omega )d\omega 
  = \frac{\mathrm{Re}(f_j)}{jK|f_j|^2}, &  \\
\end{array} \right.
\end{eqnarray*}
which determines eigenvalues of $T_j$.
In what follows, we restrict our problem to the model (\ref{KMP}), for which
the coupling function is given by $f(\theta ) = \sin (\theta +\alpha _1) + h\sin 2(\theta +\alpha _2)$.
In this case, we have
\begin{eqnarray}
f_1 = \frac{1}{2}(\sin \alpha _1 - i\cos \alpha _1) , \quad f_2 = \frac{h}{2}(\sin 2\alpha _2 - i\cos 2\alpha _2),
\label{f1f2}
\end{eqnarray}
and $f_j = 0$ for $j\neq 1, 2$.
The spectrum of the operator $T_j$ for $j\neq 1, 2$ consists of the continuous spectrum
on the imaginary axis.
The eigenvalues of $T_1$ and $T_2$ are determined by the equations
\begin{equation}
\left\{ \begin{array}{ll}
\displaystyle \int_{\R}\! \frac{x}{x^2 + (y-\omega )^2}g(\omega )d\omega
   = \frac{2}{K} \cos \alpha _1,  & \\[0.4cm]
\displaystyle \int_{\R}\! \frac{y-\omega }{x^2 + (y-\omega )^2}g(\omega )d\omega 
  = \frac{2}{K} \sin \alpha _1, &  \\
\end{array} \right.
\label{eigen1}
\end{equation}
and 
\begin{equation}
\left\{ \begin{array}{ll}
\displaystyle \int_{\R}\! \frac{x}{x^2 + (y-2\omega )^2}g(\omega )d\omega
   = \frac{h}{K} \cos 2\alpha _2,  & \\[0.4cm]
\displaystyle \int_{\R}\! \frac{y-2\omega }{x^2 + (y-2\omega )^2}g(\omega )d\omega 
  = \frac{h}{K} \sin 2\alpha _2, &  \\
\end{array} \right.
\label{eigen2}
\end{equation}
respectively.

Let us study the properties of Eq.(\ref{eigen-eq}) for $j=1$ or equivalently Eq.(\ref{eigen1}).
Define a function $D(\lambda )$ to be
\begin{equation}
D(\lambda ) = \int_{\R}\! \frac{1}{\lambda -i\omega }g(\omega )d\omega.
\label{D} 
\end{equation}
The next lemma follows from well known properties of the Poisson integral and the Hilbert transform.
See \cite{SW} for the proof.
\\[0.2cm]
\textbf{Lemma \thedef.} Suppose $g(\omega )$ is smooth at $\omega =y$. Then, the equality
\begin{eqnarray*}
\lim_{\lambda \to +0 + iy} D^{(n)}(\lambda ) 
&=& (-1)^n n! \cdot \lim_{\lambda \to +0+iy} \int_{\R}\! \frac{1}{(\lambda -i\omega )^{n+1}} g(\omega )d\omega \\
&=& \frac{1}{i^n} \cdot \lim_{\lambda \to +0+iy} \int_{\R}\! \frac{1}{\lambda -i\omega } g^{(n)}(\omega )d\omega \\
&=& \frac{1}{i^n}\left(  \pi g^{(n)} (y) - i\pi H[g^{(n)}](y)\right)
\end{eqnarray*}
holds for $n=0,1,2,\cdots $, where $\lambda \to +0+iy$ implies the limit to the point $iy \in i\R$
from the right half plane and $H[g]$ denotes the Hilbert transform defined by
\begin{eqnarray*}
H[g](y) &=& \frac{-1}{\pi} \mathrm{p.v.} \int_{\R} \frac{1}{\omega }g(\omega +y )d\omega \\
&=& \frac{-1}{\pi} \lim_{\varepsilon \to +0} \int^{\infty}_{\varepsilon }
\frac{1}{\omega }\left( g(y+\omega ) - g(y-\omega )\right) d\omega .
\end{eqnarray*}
\textbf{Lemma \thedef.} Suppose $K>0$ and $\cos \alpha _1>0$. Then,
\\
(i) If an eigenvalue $\lambda $ of $T_1$ exists, it satisfies $\mathrm{Re} (\lambda ) > 0$.
\\
(ii) If $K>0$ is sufficiently large, there exists at least one eigenvalue $\lambda $ near infinity
on the right half plane.
\\
(iii) If $K>0$ is sufficiently small, there are no eigenvalues of $T_1$.
\\[0.2cm]
\textit{Proof.} Part (i) is obvious from the first equation of (\ref{eigen1}).
For (ii), Eq.(\ref{eigen-eq}) gives $1/\lambda + O(1/\lambda ^2) = 1/iKf_1$ when $|\lambda |$ is large.
Rouch\'{e}'s theorem proves that Eq.(\ref{eigen-eq}) has a root $\lambda \sim iKf_1$ if $K>0$
is sufficiently large.
To prove (iii), we can show that the left hand side of the first equation of 
Eq.(\ref{eigen1}) is bounded for any $x, y \in \R$, while the right hand side is not as $K\to +0$.
See \cite{Chi2} for the detail. $\Box$
\\

Eq.(\ref{eigen1}) combined with Lemma 3.1 yields
\begin{eqnarray*}
\left\{ \begin{array}{ll}
\displaystyle \lim_{x\to +0} \int_{\R}\! \frac{x}{x^2 + (y-\omega )^2}g(\omega )d\omega
   = \frac{2}{K} \cos \alpha _1 = \pi g(y),  & \\[0.4cm]
\displaystyle \lim_{x\to +0} \int_{\R}\! \frac{y-\omega }{x^2 + (y-\omega )^2}g(\omega )d\omega 
  = \frac{2}{K} \sin \alpha _1 = \pi H[g](y). &  \\
\end{array} \right.
\end{eqnarray*}
Thus, we obtain $\tan \alpha _1 \cdot g(y) = H[g](y)$.
Let $y_1, y_2, \cdots $ be roots of this equation, and put $K_j = 2\cos \alpha _1/(\pi g(y_j))$.
The pair $(y_j, K_j)$ describes that some eigenvalue $\lambda =\lambda_j (K)$ of $T_1$
on the right half plane converges to the point $iy_j$ on the imaginary axis as $K \to K_j+0$.
Since $\mathrm{Re}(\lambda )>0$, the eigenvalue $\lambda_j (K)$ is absorbed into the 
continuous spectrum on the imaginary axis and disappears at $K=K_j$.
Suppose that $y_c$ satisfies $\sup_j \{ g(y_j)\} = g(y_c)$ and put
\begin{equation}
K_c = \inf_j \{ K_j\} = \frac{2\cos \alpha _1}{\pi g(y_c)}.
\end{equation}
In this section, we need not assume that $y_c$ satisfying $\sup_j \{ g(y_j)\} = g(y_c)$ is unique,
although we will assume it in Sec.6.
We have proved that 
\\[0.2cm]
\textbf{Proposition \thedef.} Suppose $g(\omega )$ is continuous and $\cos \alpha _1 > 0$.
Then, 
\\
(i) The continuous spectrum of $T_1$ is given by $\sigma _c(T_1) = i \cdot \mathrm{supp}(g) \subset i\R$.
\\
(ii) Eigenvalues of $T_1$ are given by roots of Eq.(\ref{eigen-eq}) for $j=1$.
If it exists, $\mathrm{Re}(\lambda ) > 0$.
\\
(iii) Any eigenvalue converges to some point $iy_j$ on the imaginary axis as $K\to K_j + 0 > 0$
and disappears at $K=K_j$.
\\
(iv) When $0<K<K_c$, $T_1$ has no eigenvalues on the right half plane.
\\

In what follows, $\lambda _c(K)$ denotes the eigenvalue of $T_1$ satisfying
$\lambda _c \to +0+iy_c$ as $K\to K_c+0$.
The following formulae will be used later.
\\[0.2cm]
\textbf{Lemma \thedef.} The equalities
\begin{eqnarray*}
& & D(iy_c) := \lim_{\lambda \to +0+iy_c}D(\lambda ) = \frac{1}{iK_cf_1}, \\
& & \frac{d\lambda_c }{dK}\Bigl|_{K=K_c} = \frac{-1}{iK_c^2 f_1 D'(iy_c)}
\end{eqnarray*}
hold.
\\[0.2cm]
\textit{Proof.} The first one follows from Eq.(\ref{eigen-eq}) and the definition of $(y_c, K_c)$.
The derivative of Eq.(\ref{eigen-eq}) as a function of $\lambda $ gives
\begin{eqnarray*}
D'(\lambda )= \frac{-1}{iK(\lambda )^2f_1} \frac{dK}{d\lambda }.
\end{eqnarray*}
This proves the second one. $\Box$
\\

The eigenvalues of $T_2$ are calculated in a similar manner.
The limit $x\to +0$ for Eq.(\ref{eigen2}) provides
\begin{eqnarray*}
\left\{ \begin{array}{ll}
\displaystyle \lim_{x\to +0}\int_{\R}\! \frac{x}{x^2 + (y-2\omega )^2}g(\omega )d\omega
   = \frac{h}{K} \cos 2\alpha _2 = \frac{1}{2}\pi g(y/2),  & \\[0.4cm]
\displaystyle \lim_{x\to +0}\int_{\R}\! \frac{y-2\omega }{x^2 + (y-2\omega )^2}g(\omega )d\omega 
  = \frac{h}{K} \sin 2\alpha _2 = \frac{1}{2}\pi H[g](y/2). &  \\
\end{array} \right.
\end{eqnarray*}
Let $y_1, y_2,\cdots $ be roots of the equation $\tan 2\alpha _2 \cdot g(y/2) = H[g](y/2)$.
Define $K_j^{(2)} = 2h \cos 2\alpha _2/(\pi g(y_j/2))$ and $K_c^{(2)} = \inf_j \{ K_j^{(2)}\}$.
Then, $T_2$ satisfies a similar statement to Proposition 3.3.

In what follows, we assume the following;
\\[0.2cm]
\textbf{(A1)} $\cos \alpha _1 > 0$ and $K_c < K_c^{(2)}$.
\\

The assumption $\cos \alpha _1 > 0$ shows that eigenvalues of $T_1$ lie on the right half plane (Lemma 3.2).
The assumption $K_c < K_c^{(2)}$ implies that the eigenvalue $\lambda _c$ of $T_1$ still exists
on the right half plane after all eigenvalues of $T_2$ disappear as $K$ decreases.
In other words, as $K$ increases from zero, the eigenvalue $\lambda _c$ of $T_1$ first emerges
from the imaginary axis before some eigenvalue of $T_2$ emerges.
When $\alpha _1 = 2\alpha _2$, $K_c < K_c^{(2)}$ is equivalent to the condition $h<1$
(i.e. the first harmonic $\sin \theta $ dominates the second harmonic $h\sin 2\theta $
for the attractive coupling between oscillators).
\\[0.2cm]
\textbf{Theorem \thedef\, (Instability of the incoherent state).}
\\
Suppose (A1) and $g(\omega )$ is continuous.
If $0<K<K_c$, the spectra of operators $T_1, T_2, \cdots $ consist only of the continuous spectra on the imaginary axis.
There exists a small number $\varepsilon >0$ such that
when $K_c<K<K_c+\varepsilon $, the eigenvalue $\lambda _c$ of $T_1$ lies on the right half plane.
Therefore, the incoherent state is linearly unstable.
\\

This suggests that a first bifurcation occurs at $K=K_c$ and the eigenvalue $\lambda _c$ of $T_1$
plays an important role to the bifurcation.
\\[0.2cm]
\textbf{Example \thedef.} In this example, suppose that $g$ is an even function.
Then, the Hilbert transform $H[g]$ is odd. In particular $H[g](0) = 0$.

(i) Suppose that $\alpha _1 \neq n\pi$.
Then, $y=0$ does not satisfy the equation $\tan \alpha _1 \cdot g(y) = H[g](y)$ unless $g(0) = 0$.
If $y_1 \neq 0$ satisfies this equation, $-y_1$ does not unless $g(y_1) = 0$.
Hence, a Hopf bifurcation does not occur at $K=K_c$.
We will consider this situation in this paper.

(ii) Suppose that $\alpha _1 = 0$.
Then, $y = 0$ is a root of the equation $\tan \alpha _1 \cdot g(y) = H[g](y)$.
If $y_1\neq 0$ is a root, so is $-y_1$.
If further suppose that $g(\omega )$ is a unimodal function ($g(0) > g(y_1)$),
$y_c =0$ and the eigenvalue $\lambda _c$ emerges at the origin at $K=K_c$.
As a result, a pitchfork ($h=0$) or a transcritical ($h<1$) bifurcation occurs at $K=K_c$, see Fig.\ref{fig2}.
If $g(\omega )$ is bimodal, $g(0)<g(y_c)= g(-y_c)$ in general.
As a result, a pair of eigenvalues emerges from the imaginary axis at $K=K_c$
and a Hopf bifurcation seems to occur.
In this paper, we are interested in the effect of phase-lag $\alpha _1 \neq 0$, and thus a Hopf bifurcation is not considered. 

When $g(\omega )$ is even and unimodal, the implicit function theorem shows that $y_c$
is given by
\begin{equation}
y_c = -\frac{g(0)}{H[g]'(0)}\alpha _1 + O(\alpha _1^2)
\end{equation}
as $\alpha _1 \to 0$.


\section{Linear stability}

When $0<K<K_c$, there are no spectra of operators $T_1, T_2, \cdots $ on the right half plane,
while the continuous spectra of them lie on the imaginary axis.
Hence, one may expect that the incoherent state is neutrally stable.
Nevertheless, we will show that the incoherent state is asymptotically stable in a certain weak sense.
The key idea is as follows.

For $j\geq 3$, $T_j = ij\omega $ is the multiplication operator.
Thus, the semigroup is easily obtained as $e^{T_jt} = e^{ij\omega t}$.
For any $\phi \in L^2 (\R, g(\omega )d\omega )$, $e^{T_jt}\phi$ is actually neutrally stable
with respect to the $L^2 (\R, g(\omega )d\omega )$-topology; $|| e^{T_jt}\phi || = || \phi ||$.
Let us consider the inner product with another function $\psi$
\begin{eqnarray*}
(e^{T_j t}\phi, \psi) = \int_{\R}\! e^{ij\omega t} \phi (\omega )\overline{\psi (\omega )}g(\omega )d\omega . 
\end{eqnarray*}
It is known as the Riemann-Lebesgue lemma that if the function $\phi (\omega )\overline{\psi (\omega )}g(\omega )$
has some regularity, this quantity tends to zero as $t\to \infty$.
In particular if $\phi (\omega )\overline{\psi (\omega )}g(\omega )$ is an analytic function,
$(e^{T_j t}\phi, \psi)$ decays to zero exponentially.
Thus, $e^{T_j t}\phi$ decays to zero in some weak sense.

To be precise, we need the following assumption.
Let $\delta $ be a positive number and define the stripe region on $\C$
\begin{eqnarray*}
S(\delta ) := \{ z\in \C \, | \, 0 \leq \mathrm{Im}(z) \leq \delta \}.
\end{eqnarray*}
We assume that
\\[0.2cm]
\textbf{(A2)} The density function $g(\omega )$ has an analytic continuation to the region $S(\delta )$.
On $S(\delta )$, there exists a constant $C>0$ such that the estimate
\begin{equation}
|g(z)| \leq \frac{C}{1+|z|^2}, \quad z\in S(\delta )
\end{equation}
holds.
\\

Let $H_+$ be the Hardy space:
the set of bounded holomorphic functions on the real axis and the upper half plane.
It is a subspace of $L^2 (\R, g(\omega )d\omega )$.
For $\psi \in H_+$, set $\psi^*(z) := \overline{\psi (\overline{z})}$.

A function $f_t \in  L^2 (\R, g(\omega )d\omega )$ parameterized by $t$ is said to be convergent to zero
in the weak sense if the inner product $(f_t, \psi^*)$ decays to zero as $t\to \infty$ for any $\psi \in H_+$.
Note that $P_0\in H_+$ and the order parameter is written as $\eta_1(t) = (Z_1, P_0) = (Z_1, P_0^*)$.
This means that it is sufficient to consider the stability in the weak sense for the stability of the order parameter.

For $j\geq 3$, we have
\begin{eqnarray*}
(e^{T_j t}\phi, \psi^*) = \int_{\R}\! e^{ij\omega t} \phi (\omega )\psi (\omega )g(\omega )d\omega . 
\end{eqnarray*}
(the conjugate $\psi^*$ is introduced to avoid the complex conjugate in the right hand side).
Now a standard technique of function theory shows
\begin{eqnarray*}
(e^{T_j t}\phi, \psi^*) &=& \int^{i\delta +\infty}_{i\delta -\infty}\! 
  e^{ij\omega t} \phi (\omega )\psi (\omega )g(\omega )d\omega \\
&=& e^{-j\delta t}\int_{\R}\! e^{ij\omega t} \phi (\omega +i\delta )\psi (\omega +i\delta )g(\omega +i\delta )d\omega.
\end{eqnarray*}
Hence, there is a constant $D$ such that $|(e^{T_j t}\phi, \psi^*)| < D e^{-j\delta t}$ for $t>0$
if $\phi, \psi \in H_+$ and if $g$ satisfies (A2).
Therefore, the trivial solution of the linearized system (\ref{linear}) for $j\geq 3$ is asymptotically stable
in the weak sense.

To prove a similar fact for $j=1,2$, we need the following fundamental lemma.
\\[0.2cm]
\textbf{Lemma \thedef \cite{Chi2,Chi3}.} Let $f(z)$ be a holomorphic function on the region $S(\delta )$.
Define a function $A[f](\lambda )$ of $\lambda $ to be 
\begin{eqnarray*}
A[f](\lambda ) = \int_{\R}\! \frac{1}{\lambda -i\omega }f(\omega )d\omega  
\end{eqnarray*} 
for $\mathrm{Re}(\lambda )>0$.
It has an analytic continuation $\hat{A}[f](\lambda )$ from the right half plane to the region 
$-\delta \leq \mathrm{Re}(\lambda ) \leq 0$ given by
\begin{equation}
\hat{A}[f](\lambda ) = \left\{ \begin{array}{ll}
A[f](\lambda ) & \mathrm{Re}(\lambda )>0 \\[0.2cm]
\displaystyle \lim_{\mathrm{Re}(\lambda ) \to +0}A[f](\lambda ) & \mathrm{Re}(\lambda ) = 0 \\[0.2cm]
A[f](\lambda ) + 2\pi f(-i\lambda ) & -\delta \leq \mathrm{Re}(\lambda ) < 0.
\end{array} \right.
\label{lemma4-1}
\end{equation}
\textit{Proof.} It is verified by using the formulae on the Poisson integral
\begin{eqnarray*}
& & \lim_{\lambda \to+0+iy} A[f](\lambda ) = \pi f(y) - i\pi H[f](y), \\
& & \lim_{\lambda \to-0+iy} A[f](\lambda ) = -\pi f(y) - i\pi H[f](y),
\end{eqnarray*}
see Lemma 3.1.
This shows that the first line and the third line of Eq.(\ref{lemma4-1}) coincide with one another
on the imaginary axis,
and it is continuous on the imaginary axis.
Thus, the right hand side of (\ref{lemma4-1}) is holomorphic. $\Box$
\\

\begin{figure}
\begin{center}
\includegraphics[scale=1.3]{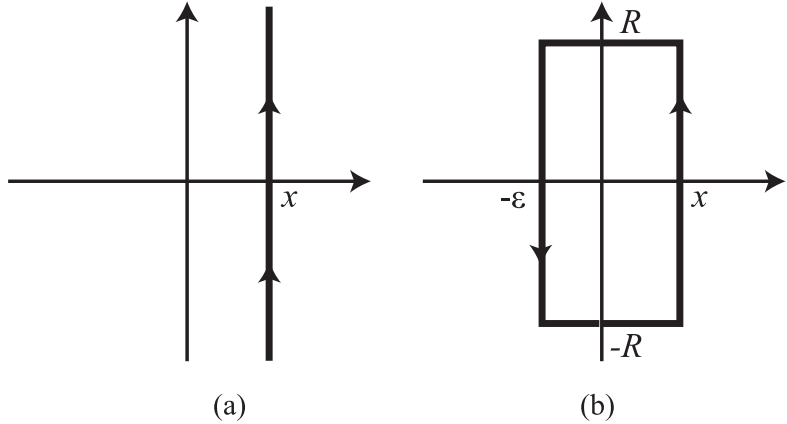}
\caption[]{Deformation of the integral path for the Laplace inversion formula.}
\label{fig3}
\end{center}
\end{figure}

It is known that the semigroup $e^{Tt}$ of an operator $T$ is expressed by the Laplace inversion formula
\begin{equation}
e^{Tt} = \lim_{y\to \infty} \frac{1}{2\pi i} \int^{x+iy}_{x-iy}\! e^{\lambda t}(\lambda -T)^{-1}d\lambda, 
\end{equation}
for $t>0$ if $T$ is a generator of a $C_0$-semigroup \cite{Yos}.
Here, $x>0$ is chosen so that the integral path is to the right of the spectrum of $T$ (see Fig.\ref{fig3}(a)).
The next purpose is to calculate the resolvent of our operator $T_1 = i\omega + iKf_1\mathcal{P}$.
\\[0.2cm]
\textbf{Lemma \thedef.} The resolvent of $T_1$ is given by
\begin{equation}
(\lambda -T_1)^{-1}\phi = (\lambda -i\omega )^{-1}\phi
 + \frac{iKf_1}{1-iKf_1D(\lambda )} ((\lambda -i\omega )^{-1}\phi, P_0) \frac{1}{\lambda -i\omega }.
\end{equation}
Let $\lambda _c$ be a simple eigenvalue of $T_1$.
The projection $\Pi_c$ to the eigenspace of $\lambda _c$ is given by
\begin{equation}
\Pi_c\phi = \frac{-1}{D'(\lambda _c)} ((\lambda _c - i\omega )^{-1} \phi, P_0) \frac{1}{\lambda _c - i\omega }.
\end{equation}
\\
\textit{Proof.} Put $R = (\lambda -T_1)^{-1}$. We have
\begin{eqnarray*}
(\lambda -T_1)R\phi = (\lambda -i\omega -iKf_1\mathcal{P})R\phi = \phi.
\end{eqnarray*}
This provides
\begin{eqnarray}
R\phi = (\lambda -i\omega )^{-1}\phi + iKf_1 (R\phi, P_0) (\lambda -i\omega )^{-1}P_0.
\label{4-6}
\end{eqnarray}
The inner product with $P_0$ yields
\begin{eqnarray*}
& & (R\phi, P_0) = ((\lambda -i\omega )^{-1}\phi, P_0) + iKf_1(R\phi, P_0) ((\lambda -i\omega )^{-1}P_0, P_0) \\
& & (R\phi, P_0) = \frac{1}{1-iKf_1 D(\lambda )}((\lambda -i\omega )^{-1}\phi, P_0).
\end{eqnarray*}
Substituting this into Eq.(\ref{4-6}) gives $R\phi$.

Next, we calculate $\Pi_c \phi$.
Let $\gamma $ be a small simple closed curve enclosing $\lambda _c$.
By the assumption, $\lambda _c$ is a pole of $(1-iKf_1 D(\lambda ))^{-1}$ of the first order,
and $(\lambda -i\omega )^{-1}$ is regular inside $\gamma $
(recall that an eigenvalue is a root of $1-iKf_1D(\lambda ) = 0$).
Hence, the Riezs projection is calculated as
\begin{eqnarray*}
\Pi_c\phi &=& \frac{1}{2\pi i} \int_{\gamma }\! (\lambda -T_1)^{-1}\phi d\lambda  \\
&=& \frac{1}{2\pi i}\cdot iKf_1 \cdot ((\lambda _c-i\omega )^{-1}\phi, P_0)\frac{1}{\lambda _c-i\omega }
\int_{\gamma }\! \frac{d\lambda }{1-iKf_1 D(\lambda )}.  
\end{eqnarray*}
The residue theorem gives
\begin{eqnarray*}
\int_{\gamma }\! \frac{d\lambda }{1-iKf_1D(\lambda )} = 2\pi i\frac{1}{-iKf_1D'(\lambda _c)}. 
\end{eqnarray*}
This completes a proof of the lemma. $\Box$
\\

Lemma 4.2 provides
\begin{eqnarray}
& & ((\lambda -T_1)^{-1}\phi, \psi^*) \nonumber \\
&=& ((\lambda -i\omega )^{-1}\phi, \psi^*)
 + \frac{iKf_1}{1-iKf_1D(\lambda )} ((\lambda -i\omega )^{-1}\phi, P_0)\cdot ((\lambda -i\omega )^{-1} \psi, P_0),
\label{resolvent}
\end{eqnarray}
which is meromorphic in $\lambda $ on the right half plane.
Suppose $\phi, \psi \in H_+$.
Due to the fundamental lemma 4.1, $((\lambda -T_1)^{-1}\phi, \psi^*)$
has an analytic continuation, possibly with new singularities, to the region 
$-\delta \leq \mathrm{Re}(\lambda ) \leq 0$
(Lemma 4.1 is applied to the factors $D(\lambda ), ((\lambda -i\omega )^{-1}\phi, \psi^*),
((\lambda -i\omega )^{-1}\phi, P_0)$ and $((\lambda -i\omega )^{-1}\psi, P_0)$).
A singularity on the left half plane is a root of the equation
\begin{equation}
1-iKf_1 (D(\lambda ) + 2\pi g(-i\lambda )) = 0.
\label{resonance}
\end{equation}
Such a singularity of the analytic continuation of the resolvent on the left half plane 
is called the generalized eigenvalue (see Sec.5 for the detail).

Now we can estimate the behavior of the semigroup by using the analytic continuation.
We have
\begin{eqnarray*}
(e^{T_1t}\phi, \psi^*) 
= \lim_{y\to \infty} \frac{1}{2\pi i} \int^{x+iy}_{x-iy}\! e^{\lambda t}((\lambda -T_1)^{-1}\phi, \psi^*) d\lambda, 
\end{eqnarray*}
where the integral path is given as in Fig.\ref{fig3} (a).
When $\phi, \psi \in H_+$, the integrand $((\lambda -T_1)^{-1}\phi, \psi^*)$ has an analytic continuation
to the region $-\delta \leq \mathrm{Re}(\lambda ) \leq 0$ which is denoted by $\mathcal{R}(\lambda )$.
\\[0.2cm]
\textbf{Lemma \thedef.}
Fix $K$ such that $0<K<K_c$.
Take positive numbers $\varepsilon, R$ and consider the rectangle shaped closed path $C$
represented in Fig.\ref{fig3} (b).
If $\varepsilon >0$ is sufficiently small, the analytic continuation of $((\lambda -T_1)^{-1}\phi, \psi^*)$
is holomorphic inside $C$ for any $R>0$.
\\

If this lemma is true, we have
\begin{eqnarray*}
0&=& \int^{x+iR}_{x-iR}\! e^{\lambda t} ((\lambda -T_1)^{-1}\phi, \psi^*) d\lambda
     + \int^{-\varepsilon -iR}_{-\varepsilon +iR} e^{\lambda t} \mathcal{R}(\lambda ) d\lambda \\
&+& \int^{iR}_{x+iR}\! e^{\lambda t} ((\lambda -T_1)^{-1}\phi, \psi^*) d\lambda
     + \int^{iR-\varepsilon }_{iR} e^{\lambda t} \mathcal{R}(\lambda ) d\lambda \\
&+& \int^{-iR+x}_{-iR}\! e^{\lambda t} ((\lambda -T_1)^{-1}\phi, \psi^*) d\lambda
     + \int^{-iR}_{-\varepsilon -iR} e^{\lambda t} \mathcal{R}(\lambda ) d\lambda .
\end{eqnarray*}
Due to the assumption (A2), we can verify that four integrals in the second and third lines above
become zero as $R\to \infty$ (see Appendix for the proof).
Thus, we obtain
\begin{eqnarray*}
(e^{T_1t}\phi, \psi^*) 
= \lim_{R\to \infty} \frac{1}{2\pi i} \int^{-\varepsilon +iR}_{-\varepsilon -iR}\! 
   e^{\lambda t}\mathcal{R}(\lambda ) d\lambda.
\end{eqnarray*}
This proves $|(e^{T_1t}\phi, \psi^*) | \sim O(e^{-\varepsilon t})$ as $t\to \infty$.
We can show the same result for the operator $T_2$.
\\[0.2cm]
\textbf{Theorem \thedef\, (Local stability of the incoherent state).}
\\
Suppose (A1) and (A2).
When $0<K<K_c$, $(e^{T_jt}\phi, \psi^*) $ decays to zero exponentially as $t\to \infty$
for any $j=1,2,\cdots $ and any $\phi, \psi\in H_+$.
Thus, the incoherent state is linearly asymptotically stable in the weak sense.
\\

\textit{Proof of Lemma 4.3.}
When $0<K<K_c$, $T_1$ has no eigenvalues on the right half plane and the imaginary axis, so that
$((\lambda -T_1)^{-1}\phi, \psi^*)$ is holomorphic on the right half plane and the imaginary axis.
A singularity on the left half plane is a root of Eq.(\ref{resonance}).
Since $D(\lambda )$ and $g(-i\lambda )$ are holomorphic,
the set of singularities has no accumulation points.
Hence, for each $R>0$, there are no singularities inside the path $C$ if $\varepsilon =\varepsilon (R)$
is sufficiently small.
Finally, Eq.(\ref{resonance}) is estimated as $1=O(1/\lambda )$ as $|\lambda | \to \infty$ in the region 
$-\delta \leq \mathrm{Re}(\lambda ) \leq 0$.
This means that there exists $R_0>0$ such that if $|\lambda | > R_0$, there are no singularities
in the region $-\delta \leq \mathrm{Re}(\lambda ) \leq 0$.
Then, the lemma holds with $\varepsilon = \varepsilon (R_0)$. $\Box$

\section{The generalized spectral theory}

For the study of a bifurcation, we need generalized spectral theory developed in \cite{Chi3}
and applied to the Kuramoto model in \cite{Chi2} because the operator $T_1$
has the continuous spectrum on the imaginary axis (thus, the standard center manifold reduction is not applicable).
In this section, a simple review of the generalized spectral theory is given.
All proofs are included in \cite{Chi2, Chi3}.
We have already encountered a part of the theory; the analytic continuation of the resolvent
and its singularity called the generalized eigenvalue.
Here, we will reformulate these concept in more functional-analytic manner.

Let $H_+$ be the Hardy space on the upper half plane with the norm
\begin{equation}
|| \phi ||_{H_+} = \sup_{\mathrm{Im}(z) > 0} |\phi (z)|.
\end{equation}
With this norm, $H_+$ is a Banach space.
Let $H_+'$ be the dual space of $H_+$; the set of continuous anti-linear functionals on $H_+$.
For $\mu \in H_+'$ and $\phi \in H_+$, $\mu (\phi )$ is denoted by $\langle \mu \,|\, \phi \rangle$.
For any $a,b \in \C,\, \phi, \psi \in H_+$ and $\mu, \xi \in H_+'$, the equalities
\begin{eqnarray*}
& & \langle \mu \,|\,  a \phi + b\psi\rangle 
   = \overline{a} \langle \mu \,|\,  \phi \rangle + \overline{b} \langle \mu \,|\, \psi \rangle, \\
& & \langle a\mu + b\xi \,|\, \phi \rangle
   = a \langle \mu \,|\, \phi \rangle + b \langle \xi \,|\, \phi \rangle,
\end{eqnarray*}
hold. An element of $H_+'$ is called a generalized function.
The space $H_+$ is a dense subspace of $L^2=L^2 (\R, g(\omega )d\omega )$
and the embedding $H_+ \hookrightarrow L^2$ is continuous.
Then, we can show that the dual $(L^2)'$ of $L^2$ is dense in $H_+'$ and it is continuously embedded in $H_+'$.
Since $L^2$ is a Hilbert space satisfying $ (L^2)' \simeq L^2$, we have
three topological vector spaces called a Gelfand triplet
\begin{eqnarray*}
H_+ \subset L^2 (\R, g(\omega )d\omega ) \subset H_+'.
\end{eqnarray*}
If an element $\phi \in H_+'$ is included in $L^2 (\R, g(\omega )d\omega )$, 
then $\langle \phi \,|\, \psi \rangle$ is given by 
\begin{eqnarray*}
\langle \phi \,|\, \psi \rangle := (\phi, \psi^*) = \int_{\R}\! \phi (\omega )\psi(\omega )g(\omega )d\omega . 
\end{eqnarray*}
Our operator $T_1$ and the above triplet satisfy all assumptions given in \cite{Chi3} to
develop a generalized spectral theory.
Now we give a brief review of the theory.
In what follows, we assume (A2).

The multiplication operator $\phi \mapsto i\omega \phi$ has the continuous spectrum 
on the imaginary axis; its resolvent is given by $(\lambda -i\omega )^{-1}$, and
it is not included in $L^2(\R, g(\omega )d\omega )$ 
when $\lambda $ is a purely imaginary number.
Nevertheless, we show that the resolvent has an analytic continuation from the right half plane to the 
left half plane in the generalized sense.
We define an operator $A(\lambda ) : H_+ \to H_+'$, parameterized by $\lambda \in \C$, to be
\begin{eqnarray*}
\langle A(\lambda )\phi \,|\, \psi \rangle
 = \left\{ \begin{array}{ll}
\displaystyle ((\lambda -i\omega )^{-1}\phi, \psi^*), & \mathrm{Re}(\lambda )>0, \\[0.4cm]
\displaystyle \displaystyle \lim_{\mathrm{Re}(\lambda ) \to +0}((\lambda -i\omega )^{-1}\phi, \psi^*)
   & \mathrm{Re}(\lambda )=0, \\[0.4cm]
\displaystyle ((\lambda -i\omega )^{-1}\phi, \psi^*) \\
\displaystyle  \quad +2\pi \phi (-i\lambda) \psi (-i\lambda)g(-i\lambda )& -\delta \leq \mathrm{Re}(\lambda )<0,
\end{array} \right.
\end{eqnarray*}
for $\phi, \psi \in H_+$.
Due to Lemma 4.1, $\langle A(\lambda )\phi \,|\, \psi \rangle$ is holomorphic.
That is, $A(\lambda )\phi$ is a $H_+'$-valued holomorphic function in $\lambda $.
In particular, $A(\lambda )$ coincides with $(\lambda -i\omega )^{-1}$ when $\mathrm{Re}(\lambda )>0$.
Since the continuous spectrum of the multiplication operator by $i\omega $ is the whole imaginary axis,
 $(\lambda -i\omega )^{-1}$ does not have an analytic continuation from the right half plane to the left half plane
as an operator on $L^2(\R, g(\omega )d\omega )$, however, 
it has a continuation $A(\lambda )$ if it is regarded as an operator from $H_+$ to $H_+'$.
$A(\lambda )$ is called the generalized resolvent of the multiplication operator by $i\omega $.

The next purpose is to define an analytic continuation of the resolvent of $T_1$ in the generalized sense.
Note that $(\lambda -T_1)^{-1}$ is rearranged as
\begin{eqnarray*}
(\lambda -i\omega -iKf_1\mathcal{P})^{-1}
= (\lambda -i\omega )^{-1}\circ (\mathrm{id} - iKf_1 \mathcal{P} (\lambda -i\omega )^{-1})^{-1}.
\end{eqnarray*}
Since the analytic continuation of $(\lambda -i\omega )^{-1}$ in the generalized sense is $A(\lambda )$,
we define the generalized resolvent $\mathcal{R}(\lambda ) : H_+ \to H_+'$ of $T_1$ by
\begin{eqnarray*}
\mathcal{R}(\lambda ):=A(\lambda ) \circ \left( \mathrm{id} - iKf_1 \mathcal{P}^\times A(\lambda ) \right)^{-1},
\end{eqnarray*}
where $\mathcal{P}^\times : H'_+ \to H'_+$ is the dual operator of $\mathcal{P}$.
For each $\phi \in H_+$, $\mathcal{R}(\lambda )\phi$ is a $H_+'$-valued meromorphic function.
It is easy to verify that when $\mathrm{Re}(\lambda )>0$, it is reduced to the usual resolvent $(\lambda -T_1)^{-1}$.
Thus, $\mathcal{R}(\lambda )$ gives a meromorphic continuation of $(\lambda -T_1)^{-1}$
from the right half plane to the left half plane as a $H_+'$-valued operator.
Again, note that $T_1$ has the continuous spectrum on the imaginary axis, so that it has no continuation
as an operator on $L^2(\R, g(\omega )d\omega )$.

A generalized eigenvalue is defined as a singularity of $\mathcal{R}(\lambda )$,
namely a singularity of $\left( \mathrm{id} - iKf_1 \mathcal{P}^\times A(\lambda ) \right)^{-1}$.
\\[0.2cm]
\textbf{Definition \thedef.}
If the equation
\begin{equation}
(\mathrm{id} - iKf_1 \mathcal{P}^\times A(\lambda ))\mu = 0, \quad -\delta \leq \mathrm{Re}(\lambda )
\end{equation}
has a nonzero solution $\mu$ in $H_+'$ for some $\lambda \in \C$, $\lambda $ is called a generalized eigenvalue
and $\mu$ is called a generalized eigenfunction.
\\

It is easy to verify that this equation is equivalent to
\begin{eqnarray}
\frac{1}{iKf_1} = \left\{ \begin{array}{ll}
D(\lambda ) & \mathrm{Re}(\lambda )>0,  \\
\displaystyle \lim_{\mathrm{Re}(\lambda ) \to +0}D(\lambda ) & \mathrm{Re}(\lambda )=0,  \\
D(\lambda )+2\pi g(-i\lambda ) & -\delta \leq \mathrm{Re}(\lambda )<0.
\end{array} \right.
\label{eigeneq3}
\end{eqnarray}
When $\mathrm{Re}(\lambda ) > 0$, this is reduced to Eq.(\ref{eigen-eq}) with $j=1$.
In this case, $\mu$ is included in $L^2(\R, g(\omega )d\omega )$ and 
a generalized eigenvalue on the right half plane is an eigenvalue in the usual sense.
When $\mathrm{Re}(\lambda ) \leq 0$, this equation is equivalent to Eq.(\ref{resonance}).
The associated generalized eigenfunction is not included in $L^2 (\R, g(\omega )d\omega )$
but an element of the dual space $H'_+$.
Although a generalized eigenvalue is not a true eigenvalue of $T_1$, 
it is an eigenvalue of the dual operator:
\\[0.2cm]
\textbf{Theorem \thedef\, \cite{Chi2,Chi3}.}
Let $\lambda $ and $\mu$ be a generalized eigenvalue and the associated generalized eigenfunction.
The equality $T_1^\times \mu = \lambda \mu$ holds.
\\ 

Let $\lambda _0$ be a generalized eigenvalue of $T_1$ and 
$\gamma _0$ a small simple closed curve enclosing $\lambda _0$.
The generalized Riesz projection $\Pi_0 : H_+\to H_+'$ is defined by
\begin{eqnarray*}
\Pi_0 = \frac{1}{2\pi i}\int_{\gamma _0}\! \mathcal{R}(\lambda ) d\lambda .
\end{eqnarray*}
As in the usual spectral theory, the image of it gives the generalized eigenspace associated with $\lambda _0$,
and it satisfies $\Pi_0 T_1^\times = T_1^\times \Pi_0$ \cite{Chi3}.

\begin{figure}
\begin{center}
\includegraphics[scale=1.3]{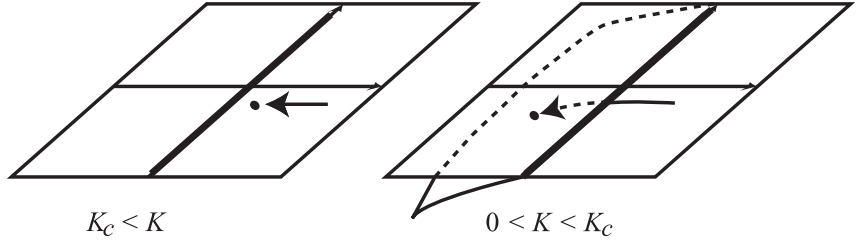}
\caption[]{The motion of the (generalized) eigenvalue as $K$ decreases.
When $0<K<K_c$, it lies on the second Riemann sheet of the resolvent and it is not a usual eigenvalue
but a generalized eigenvalue.}
\label{fig4}
\end{center}
\end{figure}

Let $\lambda = \lambda _c(K)$ be an eigenvalue of $T_1$ defined in Sec.3.
Recall that when $K_c<K$, $\lambda _c$ exists on the right half plane.
As $K$ decreases, $\lambda _c$ goes to the left side, and
at $K= K_c$, $\lambda _c$ is absorbed into the continuous spectrum on the imaginary axis and disappears.
However, we can show that even for $0<K<K_c$, 
$\lambda _c$ remains to exist as a root of Eq.(\ref{eigeneq3}) because the right hand side of Eq.(\ref{eigeneq3})
is holomorphic.
This means that although $\lambda _c$ disappears from the original complex plane at $K=K_c$,
it still exists for $0<K<K_c$ as a generalized eigenvalue on the Riemann surface of the generalized resolvent 
$\mathcal{R}(\lambda )$.
In the generalized spectral theory, the resolvent $(\lambda -T_1)^{-1}$ is regarded as an operator 
from $H_+$ to $H_+'$, not on $L^2(\R , g(\omega )d\omega )$.
Then, it has an analytic continuation from the right half plane to the left half plane
as $H_+'$-valued operator.
The continuous spectrum on the imaginary axis becomes a branch cut of the Riemann surface of the resolvent.
On the Riemann surface, the left half plane is two-sheeted (see Fig.\ref{fig4}).
We call a singularity of the generalized resolvent on the second Riemann sheet the generalized eigenvalue. 

On the dual space $H_+'$, the weak dual topology (weak star topology) is equipped;
a sequence $\{ \mu_n \} \subset H_+'$ is said to be convergent to $\mu \in H_+'$
if $\langle \mu_n \,|\, \psi \rangle \in \C$ is convergent to $\langle \mu \,|\, \psi \rangle$ for each $\psi \in H_+$.
Recall that an eigenfunction of a usual eigenvalue $\lambda $ of $T_1$ is given by 
$v_\lambda (\omega ) = (\lambda -i\omega )^{-1}$ (Eq.(\ref{ef})).
A generalized eigenfunction $\mu_\lambda $ of a generalized eigenvalue $iy$ on the imaginary axis is given by
\begin{eqnarray*}
\mu_\lambda = \lim_{\lambda \to +0+iy} \frac{1}{\lambda -i\omega },
\end{eqnarray*}
where the limit is considered with respect to the weak dual topology (compare with Eq.(\ref{ef})).
This means that $\langle \mu_\lambda  \,|\, \psi \rangle$ is defined by
\begin{equation}
\langle \mu_\lambda  \,|\, \psi \rangle
 = \lim_{\lambda \to +0+iy} \langle \frac{1}{\lambda -i\omega } \,|\, \psi \rangle
 = \lim_{\lambda \to +0+iy} \int_{\R}\! \frac{1}{\lambda -i\omega }\psi (\omega )g(\omega )d\omega . 
\end{equation}
A generalized eigenfunction $\mu_\lambda $ associated with a generalized eigenvalue $\lambda $
on the left half plane is given by
\begin{equation}
\langle \mu_\lambda  \,|\, \psi \rangle
 =  \int_{\R}\! \frac{1}{\lambda -i\omega }\psi (\omega )g(\omega )d\omega + 2\pi \psi (-i\lambda )g(-i\lambda ). 
\end{equation}

To perform a center manifold reduction, we need the definition of a center subspace.
Usually, it is defined to be an eigenspace associated with eigenvalues on the imaginary axis.
For our case, the operators $T_1, T_2,\cdots $ have the continuous spectra on the imaginary axis.
Therefore, we define a generalized center subspace as a space spanned by generalized eigenfunctions
associated with generalized eigenvalues on the imaginary axis.
Note that this is a subspace of the dual $H_+'$, not of $L^2(\R, g(\omega )d\omega )$.
As $K$ increases from zero, one of the generalized eigenvalues of $T_1$ gets across the imaginary axis
at $K=K_c$, and it becomes a usual eigenvalue on the right half plane (see Fig.\ref{fig4}).
Hence, there is a nontrivial generalized center subspace at $K=K_c$ given by
\begin{eqnarray*}
\mathbf{E}^c := \mathrm{span}\{ \mu_\lambda \, | \, \lambda _c = iy_c \},
\end{eqnarray*}
where $\lambda _c = iy_c$ is defined as in Section 3.
The next goal is to perform a center manifold reduction.

In Chiba \cite{Chi2}, the existence of the one dimensional center manifold in $H'_+$
that is tangent to $\mathbf{E}^c$ is proved
for the standard Kuramoto model ($\alpha _1 = h=0$).
By restricting the continuous model onto the center manifold, the dynamics (\ref{1-4})
of the order parameter is derived.
In this paper, we formally perform the center manifold reduction without a proof of the 
existence of a center manifold.


\section{Center manifold reduction}

Recall that $y_c \in \R$ is defined as a number satisfying $\sup_j \{ g(y_j)\} = g(y_c)$,
where $y_1, y_2,\cdots $ are roots of $\tan \alpha _1 \cdot g(y) = H[g](y)$.
This gives a point $iy_c$ on the imaginary axis that the eigenvalue $\lambda = \lambda _c(K)$
of $T_1$ approaches as $K\to K_c+0$.
For a bifurcation, we assume the following:
\\[0.2cm]
\textbf{(A3)} $y_c$ is uniquely determined.
The corresponding eigenvalue $\lambda _c(K)$ of $T_1$ is simple near $K_c$.
\\[0.2cm]
\textbf{(A4)} The real part of $\displaystyle \frac{d\lambda_c}{dK}\Bigl|_{K=K_c}$ is positive.
\\

The assumption (A3) implies that the generalized center subspace at $K=K_c$ is a one dimensional space given by 
\begin{eqnarray*}
\mathbf{E}^c = \mathrm{span}\{ \mu_c := \lim_{\lambda \to +0+iy_c}\frac{1}{\lambda -i\omega } \}.
\end{eqnarray*}
The assumption (A4) means that the generalized eigenvalue $\lambda _c$ of $T_1$ transversely 
gets across the imaginary axis from the left to the right.
\\[0.2cm]
\textbf{Remark.} (A3) and (A4) are assumptions for $g(\omega )$.
If $g$ is even and unimodal (unimodal means that $g(y_1) < g(y_2)$ when $y_1 < y_2 < 0$, and 
$g(y_1) > g(y_2)$ when $0<y_1 <y_2$), (A3) is automatically satisfied.
\\

In what follows, we assume (A1) to (A4).
We expect that a one dimensional bifurcation occurs at $K=K_c$.
Eventually we will find that oscillators approximately rotate with the velocity $y_c$ (Theorem 1.5).
Hence, we introduce the moving coordinates by $\widetilde{\theta } = \theta -y_ct$ 
and $\widetilde{\omega } = \omega -y_c$.
Put
\begin{equation}
\widetilde{g}(\widetilde{\omega }):= g(\widetilde{\omega }+y_c) = g(\omega ),
\quad \zeta_l := e^{-ily_ct}\eta_l, \quad Y_l := e^{-ily_ct}Z_l.
\end{equation}
Then, Eq.(\ref{Z}) is rewritten as
\begin{equation}
\dot{Y}_j = \widetilde{T}_jY_j + ijK\sum_{l\neq j} f_l\zeta_lY_{j-l},
\label{Y}
\end{equation}
where $\widetilde{T}_j$ is defined by replacing $g(\omega )$ in $T_j$ with $\widetilde{g}(\widetilde{\omega })$
(i.e. $g(\omega )$ in the definition of the integral operator $\mathcal{P}$ 
is replaced by $\widetilde{g}(\widetilde{\omega })$).
It is easy to verify that if $y_c$ is a root of the equation $\tan \alpha _1 \cdot g(y) = H[g](y)$,
then $y=0$ satisfies the equation $\tan \alpha _1 \cdot \widetilde{g}(y) = H[\widetilde{g}](y)$.
Therefore, the eigenvalue $\lambda _c(K)$ of the new operator $\widetilde{T}_1$
converges to the origin as $K\to K_c + 0$.

In what follows, put $\varepsilon = K-K_c$, which is a bifurcation parameter.
Our ingredients are;
\\[0.2cm]
\textbf{Equations:} The equations (\ref{Y}) for $j=1,2$ are given by
\begin{equation}
\left\{ \begin{array}{l}
\dot{Y}_1 = \widetilde{T}_cY_1 + i\varepsilon f_1\mathcal{P}Y_1 
 + iK (f_2 \zeta_2 Y_{-1} + \overline{f_1} \overline{\zeta_1} Y_2 + \overline{f_2} \overline{\zeta_2} Y_3), \\[0.2cm]
\dot{Y}_2 = \widetilde{T}_2Y_2 + 2iK (f_1\zeta_1Y_1 + \overline{f_1} \overline{\zeta_1} Y_3 + \overline{f_2} \overline{\zeta_2} Y_4),
\end{array} \right.
\label{Y2}
\end{equation}
where $\widetilde{T}_c$ is an operator $\widetilde{T}_1$ estimated at $K = K_c$ 
(i.e. $K$ is denoted as $K = K_c + \varepsilon $ and 
$\widetilde{T}_1 = \widetilde{T}_c + i\varepsilon f_1\mathcal{P}$).
Eventually we will find that the equations for $Y_3, Y_4, \cdots $ do not affect a bifurcation at $K=K_c$.
$f_1$ and $f_2$ are given by Eq.(\ref{f1f2}).
\\[0.2cm]
\textbf{Center subspace:} As $K$ increases from zero, one of the generalized eigenvalues of 
$\widetilde{T}_1$ denoted by $\lambda _c (K)$ gets across the imaginary axis at $0$ when $K=K_c$,
and it becomes a usual eigenvalue on the right half plane when $K>K_c$.
The associated generalized eigenfunction at $K=K_c$ is given by 
\begin{eqnarray*}
\mu_c = \lim_{\lambda \to +0} \frac{1}{\lambda -i\widetilde{\omega} },
\end{eqnarray*}
where the limit is taken with respect to the weal dual topology.
The generalized center subspace is spanned by $\mu_c$ in $H_+'$.
\\[0.2cm]
\textbf{Projection:} The projection to an eigenspace is given in Lemma 4.2.
The projection to the generalized center subspace spanned by $\mu_c$ at $K=K_c$ is
\begin{equation}
\Pi_c \phi = \frac{-1}{D'(iy_c)} \lim_{\lambda _c \to +0} 
 ((\lambda _c - i\widetilde{\omega } )^{-1} \phi, P_0) \mu_c.
\label{6-4}
\end{equation}

We divide our result into two cases, $h=0$ and $h\neq 0$ because 
types of bifurcations of them are different; when $h=0$, it is a pitchfork bifurcation,
while when $h\neq 0$, it is a transcritical bifurcation.

\subsection{Center manifold reduction $(h=0)$}

Suppose $h= f_2 = 0$.
We put $\Pi_cY_1 = K_c\alpha (t)\mu_c/2$.
The scalar valued function $\alpha (t)$ denotes a coordinate on the center subspace, 
and our goal is to derive the dynamics of $\alpha $.
Since the center manifold is tangent to $\mathbf{E}^c$, we write
\begin{equation}
Y_1 = \frac{K_c}{2}\alpha (t)\mu_c + O(\alpha ^2), \quad Y_j \sim O(\alpha ^2) ,\,\, (j\geq 2).
\label{6-5}
\end{equation}
We make the following ansatz
\begin{eqnarray*}
\varepsilon \sim O(\alpha ^2), \quad \frac{d\alpha }{dt} \sim O(\alpha ^2),
\end{eqnarray*}
which will be confirmed if the dynamics of $\alpha $ is derived.
Then, $\zeta_1$ is given by
\begin{eqnarray}
\zeta_1 &=& \int_{\R}\! Y_1 \cdot \widetilde{g}(\widetilde{\omega })d\widetilde{\omega }
 = \frac{K_c}{2}\alpha \lim_{\lambda \to +0} \int_{\R}\! 
     \frac{1}{\lambda -i\widetilde{\omega }}\widetilde{g}(\widetilde{\omega })d\widetilde{\omega }
      + O(\alpha ^2) \nonumber \\
&=& \frac{1}{2if_1}\alpha  + O(\alpha ^2),
\label{6-6}
\end{eqnarray}
where we have used Lemma 3.4;
\begin{eqnarray*}
\lim_{\lambda \to +0} \int_{\R}\! \frac{1}{\lambda -i\widetilde{\omega }} \widetilde{g}(\widetilde{\omega })d\widetilde{\omega }
= \lim_{\lambda \to +0+iy_c} \int_{\R}\! \frac{1}{\lambda -i\omega }g(\omega )d\omega  = \frac{1}{iK_cf_1}.
\end{eqnarray*} 
Since $dY_2/dt \sim O(\alpha ^3)$ by the ansatz, substituting (\ref{6-5}), (\ref{6-6})
into the second equation of (\ref{Y2}) yields
\begin{eqnarray}
0 = \widetilde{T}_2 Y_2 + \frac{K_c^2}{2} \alpha ^2 \mu_c + O(\alpha ^3).
\label{6-7}
\end{eqnarray}
Since $h=0$, $\widetilde{T}_2$ is the multiplication operator given by $\widetilde{T}_2 = 2i \widetilde{\omega} $.
Its inverse is given as follows.
\\[0.2cm]
\textbf{Lemma \thedef.}
Define the operator $\widetilde{T}_2^{-1} : H_+\to H_+'$ to be
\begin{eqnarray*}
\langle \widetilde{T}_2^{-1}\phi \,|\, \psi^* \rangle
 = -\frac{1}{2}\lim_{\lambda \to +0} \int\! \frac{1}{\lambda -i\widetilde{\omega} }
    \phi (\widetilde{\omega} )\psi (\widetilde{\omega }) \widetilde{g}(\widetilde{\omega })d\widetilde{\omega}.
\end{eqnarray*}
Then, it satisfies $\widetilde{T}_2^{-1} \circ \widetilde{T}_2\phi = \phi$ for $\phi \in H_+$.

See \cite{Chi2} for the proof.
It is written in operator form as
\begin{eqnarray*}
\widetilde{T}_2^{-1} =-\frac{1}{2}\lim_{\lambda \to +0}\frac{1}{\lambda -i\widetilde{\omega} }.
\end{eqnarray*}
This is applied to (\ref{6-7}) to yield
\begin{eqnarray*}
Y_2 = \frac{K_c^2}{4}\alpha ^2 \lim_{\lambda \to +0}\frac{1}{\lambda -i\widetilde{\omega} }\mu_c + O(\alpha ^3)
 = \frac{K_c^2}{4}\alpha ^2 \lim_{\lambda \to +0}\frac{1}{(\lambda -i\widetilde{\omega})^2 } + O(\alpha ^3).
\end{eqnarray*}
The projection $\Pi_c$ given in (\ref{6-4}) is applied to give
\begin{eqnarray*}
\Pi_cY_2 &=& \frac{K_c^2}{4}\alpha ^2\cdot \frac{-1}{D'(iy_c)} \cdot \lim_{\lambda \to 0}
   \int_{\R}\! \frac{1}{(\lambda -i\widetilde{\omega })^3} 
      \widetilde{g}(\widetilde{\omega })d\widetilde{\omega} \cdot \mu_c + O(\alpha ^3) \\
&=& \frac{-K_c^2}{4D'(iy_c)}\alpha ^2  
       \lim_{\lambda \to +0+iy_c} \int_{\R}\! \frac{1}{(\lambda -i\omega )^3} g(\omega )d\omega \cdot\mu_c + O(\alpha ^3) \\
&=& \frac{-K_c^2 D''(iy_c)}{8D'(iy_c)}\alpha ^2 \cdot \mu_c + O(\alpha ^3).
\end{eqnarray*}
Next, the projection is applied to the first equation of (\ref{Y2}) with $f_2 = 0$.
Since $\widetilde{T}_c^\times \Pi_c = \Pi_c \widetilde{T}_c^\times $ holds \cite{Chi3},
\begin{eqnarray*}
\frac{K_c}{2}\dot{\alpha }\mu_c
&=& \widetilde{T}_c^\times \Pi_c Y_1 + i\varepsilon f_1 \cdot \zeta_1 \cdot \Pi_c P_0 + iK\overline{f_1 \zeta_1} \Pi_cY_2 \\
&=& \widetilde{T}_c^\times \left( \frac{K_c}{2}\alpha \mu_c \right) 
    + \frac{\varepsilon \alpha }{2} \cdot \frac{-1}{D'(iy_c)} \lim_{\lambda \to +0} \int_{\R}\! 
      \frac{1}{\lambda -i\widetilde{\omega }} \widetilde{g}(\widetilde{\omega })d\widetilde{\omega } \cdot \mu_c \\
& & \quad -\frac{K_c}{2}\overline{\alpha }\cdot \frac{-K_c^2 D''(iy_c)}{8D'(iy_c)}\alpha ^2 \mu_c + O(\alpha ^4).
\end{eqnarray*}
Theorem 5.2 shows $\widetilde{T}_c^\times \mu_c = 0$.
Therefore, we obtain the dynamics on the center manifold as
\begin{eqnarray}
\frac{d\alpha }{dt} &=& \frac{-1}{iK_c^2f_1 D'(iy_c)}\varepsilon \alpha 
  + \frac{K_c^2 D''(iy_c)}{8D'(iy_c)}\alpha |\alpha |^2 + O(\alpha ^4).
\end{eqnarray}
Define constants $p_1$ and $p_3$ by
\begin{eqnarray*}
p_1 = \frac{-1}{iK_c^2f_1 D'(iy_c)}, \quad p_3=\frac{K_c^2 D''(iy_c)}{8D'(iy_c)}.
\end{eqnarray*}
Since $\varepsilon =K-K_c \sim O(\alpha ^2)$, putting $\alpha =re^{i\psi}$ yields
\begin{equation}
\left\{ \begin{array}{l}
\displaystyle \dot{r}= \mathrm{Re}(p_1)r \left( K-K_c + \frac{\mathrm{Re}(p_3)}{\mathrm{Re}(p_1)}r^2 \right)
  + O(r^4)  \\
\dot{\psi} \sim O(r^2) = O(K-K_c).  \\
\end{array} \right.
\end{equation}
This proves that the equation of $r$ has the fixed point
\begin{equation}
r_0 = \sqrt{-\frac{\mathrm{Re}(p_1)}{\mathrm{Re}(p_3)}} \cdot \sqrt{K-K_c} + O(K-K_c).
\end{equation}
Lemma 3.4 with the assumption (A4) implies $\mathrm{Re}(p_1)>0$.
Thus, if $\mathrm{Re}(p_3) < 0$, the bifurcation is supercritical; 
the fixed point exists for $K>K_c$ and is stable.
If $\mathrm{Re}(p_3) > 0$, the bifurcation is subcritical; 
the fixed point exists for $K<K_c$ and is unstable.

\subsection{Center manifold reduction $(h\neq 0)$}

Suppose $h\neq 0$.
As before, we assume Eq.(\ref{6-5}).
Then, Eqs.(\ref{6-6}), (\ref{6-7}) again hold.
The inverse of  $\widetilde{T}_2$ is given by
\\[0.2cm]
\textbf{Lemma \thedef.}
Define the operator $\widetilde{T}_2^{-1} : H_+\to H_+'$ to be
\begin{eqnarray*}
\widetilde{T}_2^{-1}\phi
 &=& -\frac{1}{2}\lim_{\lambda \to +0} \frac{1}{\lambda -i\widetilde{\omega} }\cdot \phi (\widetilde{\omega }) \\
 & & - \frac{1}{2}\lim_{\lambda \to +0}
        \frac{iKf_2}{1-iKf_2 D(iy_c)} \int_{\R}\! \frac{1}{\lambda -i\widetilde{\omega }}\phi (\widetilde{\omega })
 \widetilde{g}(\widetilde{\omega })d\widetilde{\omega }\cdot \frac{1}{\lambda -i\widetilde{\omega }}.   
\end{eqnarray*}
Then, it satisfies $\widetilde{T}_2^{-1} \circ \widetilde{T}_2\phi = \phi$ for $\phi \in H_+$.

It is verified by a straightforward calculation.
Using it, we have
\begin{eqnarray*}
Y_2 &=& \frac{K_c^2}{4}\alpha ^2 \lim_{\lambda \to +0}\frac{1}{(\lambda -i\widetilde{\omega})^2 } \\
& & + \frac{K_c^2}{4}\alpha ^2 \lim_{\lambda \to +0} \frac{iKf_2}{1-iKf_2 D(iy_c)} 
   \int_{\R}\! \frac{1}{(\lambda -i\widetilde{\omega })^2}
 \widetilde{g}(\widetilde{\omega })d\widetilde{\omega }\cdot \frac{1}{\lambda -i\widetilde{\omega }}  + O(\alpha ^3).
\end{eqnarray*}
Thus, $\zeta_2 = (Y_2, P_0)$ is calculated as
\begin{eqnarray*}
\zeta_2 &=& \frac{K_c^2}{4}\alpha ^2 \lim_{\lambda \to +0} \int_{\R}\! \frac{1}{(\lambda -i\widetilde{\omega})^2 }
            \widetilde{g}(\widetilde{\omega })d\widetilde{\omega } \\
&+ & \frac{K_c^2}{4}\alpha ^2 \lim_{\lambda \to +0} \frac{iKf_2}{1-iKf_2 D(iy_c)} 
   \int_{\R}\! \frac{1}{(\lambda -i\widetilde{\omega })^2}
 \widetilde{g}(\widetilde{\omega })d\widetilde{\omega }
 \int_{\R}\! \frac{1}{\lambda -i\widetilde{\omega }} \widetilde{g}(\widetilde{\omega })d\widetilde{\omega }+ O(\alpha ^3) \\
&=& \frac{K_c^2}{4}\frac{-D'(iy_c)}{1-iK_c f_2D(iy_c)} \alpha ^2 + O(\alpha ^3).
\end{eqnarray*}
Next, the projection is applied to the first equation of (\ref{Y2}).
\begin{eqnarray*}
\frac{K_c}{2}\dot{\alpha }\mu_c
&=& i\varepsilon f_1 \cdot \zeta_1 \cdot \Pi_c P_0 + iK_cf_2\cdot \zeta_2 \cdot \Pi_c Y_{-1} + O(\alpha ^3) \\
&=& \frac{-1}{2iK_c f_1D'(iy_c)} \varepsilon \alpha \cdot \mu_c - \frac{iK_c^3 f_2}{4}
\frac{D'(iy_c)}{1-f_2/f_1}\alpha ^2 \Pi_cY_{-1} + O(\alpha ^3).
\end{eqnarray*}
\\[0.2cm]
\textbf{Lemma \thedef.} 
$\Pi_cY_{-1}$ is given by
\begin{equation}
\Pi_cY_{-1} = \frac{-4 \cos \alpha _1}{K_c D'(iy_c)}e^{-i\mathrm{arg}(\alpha )}\mu_c + O(\alpha ).
\end{equation}
See \cite{Chi1} for the proof.
Therefore, we obtain the dynamics on the center manifold as
\begin{eqnarray}
\frac{d\alpha }{dt} &=& \frac{-1}{iK_c^2f_1 D'(iy_c)}\varepsilon \alpha 
  +\frac{2iK_cf_2 \cos \alpha _1}{1-f_2/f_1}\alpha ^2 e^{-i\mathrm{arg}(\alpha )} + O(\alpha ^3).
\end{eqnarray}
Define constants $p_1$ and $p_2$ by
\begin{eqnarray*}
p_1 = \frac{-1}{iK_c^2f_1 D'(iy_c)}, \quad p_2 = \frac{2iK_cf_2 \cos \alpha _1}{1-f_2/f_1}.
\end{eqnarray*}
Putting $\alpha =re^{i\psi}$ yields
\begin{equation}
\left\{ \begin{array}{l}
\displaystyle \dot{r}= \mathrm{Re}(p_1)r \left( K-K_c + \frac{\mathrm{Re}(p_2)}{\mathrm{Re}(p_1)}r \right)+ O(r^3)  \\
\dot{\psi} \sim O(r) = O(K-K_c).  \\
\end{array} \right.
\end{equation}
This proves that the equation of $r$ has the fixed point
\begin{equation}
r_0 = -\frac{\mathrm{Re}(p_1)}{\mathrm{Re}(p_2)} \cdot (K-K_c) + O((K-K_c)^2).
\end{equation}
Lemma 3.4 with the assumption (A4) implies $\mathrm{Re}(p_1)>0$.
Thus, if $\mathrm{Re}(p_2) < 0$, the bifurcation is supercritical; 
the fixed point exists for $K>K_c$ and is stable.
If $\mathrm{Re}(p_2) > 0$, the bifurcation is subcritical; 
the fixed point exists for $K<K_c$ and is unstable.
\\

For both cases ($h=0$ and $h\neq 0$), the order parameter is written as
\begin{eqnarray*}
\eta_1(t) &=& e^{iy_ct}\zeta_1 = e^{iy_ct}\left( \frac{1}{2if_1}\alpha +O(\alpha ^2) \right) \\
&=& \frac{r}{2if_1}e^{i(y_ct + \psi)} + O(r^2) \\
&=& re^{i\alpha _1} \cdot e^{i(y_c + O(K-K_c))t} + O(r^2).
\end{eqnarray*}
This gives Theorem 1.5.


\appendix

\section{Proof of Theorem 4.4.}

The semigroup generated by the operator $T_1$ satisfies
\begin{eqnarray*}
(e^{T_1t}\phi, \psi^*) 
= \lim_{y\to \infty} \frac{1}{2\pi i} \int^{x+iy}_{x-iy}\! e^{\lambda t}((\lambda -T_1)^{-1}\phi, \psi^*) d\lambda, 
\end{eqnarray*}
where the integral path is given as in Fig.\ref{fig3} (a).
To complete the proof of Theorem 4.4, we show that this quantity decays to zero exponentially as $t\to \infty$
when $\phi, \psi \in H_+$.

Due to Lemma 4.1, the integrand $((\lambda -T_1)^{-1}\phi, \psi^*)$ has a meromorphic continuation
to the region $-\delta \leq \mathrm{Re}(\lambda ) \leq 0$ which is denoted by $\mathcal{R}(\lambda )$.
Lemma 4.3 implies that $\mathcal{R}(\lambda )$ is holomorphic inside the closed curve $C$
represented in Fig.\ref{fig3} (b) for any $R>0$ if  $\varepsilon >0$ is sufficiently small.
Hence, Cauchy's integral theorem provides
\begin{eqnarray}
0&=& \int^{x+iR}_{x-iR}\! e^{\lambda t} ((\lambda -T_1)^{-1}\phi, \psi^*) d\lambda
     + \int^{-\varepsilon -iR}_{-\varepsilon +iR} e^{\lambda t} \mathcal{R}(\lambda ) d\lambda \nonumber \\
&+& \int^{iR}_{x+iR}\! e^{\lambda t} ((\lambda -T_1)^{-1}\phi, \psi^*) d\lambda
     + \int^{iR-\varepsilon }_{iR} e^{\lambda t} \mathcal{R}(\lambda ) d\lambda \nonumber \\
&+& \int^{-iR+x}_{-iR}\! e^{\lambda t} ((\lambda -T_1)^{-1}\phi, \psi^*) d\lambda
     + \int^{-iR}_{-\varepsilon -iR} e^{\lambda t} \mathcal{R}(\lambda ) d\lambda .
\label{A0}
\end{eqnarray}
If four integrals in the second and third lines above become zero as $R\to \infty$, we obtain
\begin{eqnarray*}
(e^{T_1t}\phi, \psi^*) 
&=& \lim_{R\to \infty} \frac{1}{2\pi i} \int^{-\varepsilon +iR}_{-\varepsilon -iR}\! 
   e^{\lambda t}\mathcal{R}(\lambda ) d\lambda \\
&=&  \lim_{R\to \infty} \frac{e^{-\varepsilon t}}{2\pi i} \int^{iR}_{-iR}\! 
   e^{\lambda t}\mathcal{R}(\lambda -\varepsilon ) d\lambda.
\end{eqnarray*}
This proves $|(e^{T_1t}\phi, \psi^*) | \sim O(e^{-\varepsilon t})$ as $t\to \infty$.

Let us show that the integral $\int^{iR-\varepsilon }_{iR} e^{\lambda t} \mathcal{R}(\lambda ) d\lambda$
tends to zero as $R\to \infty$.
It is written as
\begin{eqnarray*}
\int^{iR-\varepsilon }_{iR} e^{\lambda t} \mathcal{R}(\lambda ) d\lambda
 &=& e^{iRt} \int^{-\varepsilon }_{0}\! e^{\lambda t} \mathcal{R}(\lambda +iR)d\lambda \\
&=& e^{iRt} \int^{-\varepsilon }_{-\varepsilon _0}\! e^{\lambda t} \mathcal{R}(\lambda +iR)d\lambda
 + e^{iRt} \int^{-\varepsilon_0 }_{0}\! e^{\lambda t} \mathcal{R}(\lambda +iR)d\lambda,
\end{eqnarray*}
where $\varepsilon _0$ is an arbitrarily small number such that $0<\varepsilon _0<\varepsilon $.
Since $\mathcal{R}(\lambda )$ is holomorphic around the imaginary axis, 
there exists $M_1> 0$ such that $|e^{\lambda t} \mathcal{R}(\lambda +iR)| < M_1$ for $-\varepsilon _0 < \lambda <0$.
This shows 
\begin{equation}
\left| e^{iRt} \int^{-\varepsilon_0 }_{0}\! e^{\lambda t} \mathcal{R}(\lambda +iR)d\lambda \right| < M_1 \varepsilon _0.
\label{A1}
\end{equation}
Thus, it is sufficient to prove $\mathcal{R}(\lambda +iR) \to 0$ as $R\to \infty$ uniformly in 
$-\varepsilon < \lambda < -\varepsilon _0$.
By applying Lemma 4.1 to $((\lambda -T_1)^{-1}\phi, \psi^*)$ shown in (\ref{resolvent}),
it turns out that $\mathcal{R}(\lambda )$ is given by
\begin{eqnarray*}
& & \mathcal{R}(\lambda ) = ((\lambda -i\omega )^{-1}\phi, \psi^*) + 2\pi \phi (-i\lambda )\psi (-i\lambda )g(-i\lambda ) \\
&+& \frac{iKf_1}{1-iKf_1D(\lambda ) - 2\pi i K f_1 g(-i\lambda )} \times \\
& & \Bigl( ((\lambda -i\omega )^{-1}\phi, P_0) + 2\pi \phi (-i\lambda )g(-i\lambda ) \Bigr) \cdot
\Bigl( ((\lambda -i\omega )^{-1}\psi, P_0) + 2\pi \psi (-i\lambda )g(-i\lambda ) \Bigr).
\end{eqnarray*}
The first term $((\lambda -i\omega )^{-1}\phi, \psi^*)$ with $\lambda \mapsto \lambda +iR$ is given by
\begin{eqnarray*}
((\lambda+iR -i\omega )^{-1}\phi, \psi^*)
= \int_{\R}\! \frac{1}{\lambda -i(\omega -R)} \phi(\omega )\psi (\omega )g(\omega )d\omega.
\end{eqnarray*}
Because of the assumption (A2), there exists $C_1 > 0$ such that it is estimated as
\begin{eqnarray*}
|((\lambda+iR -i\omega )^{-1}\phi, \psi^*)|
&\leq & \int_{\R}\!\frac{C_1}{|\lambda - i(\omega -R)|} \cdot \frac{1}{1+|\omega |^2}d\omega 
\end{eqnarray*}
Since the integral in the right hand side exists, for any $\varepsilon _0>0$,
there exists $L > 0$ such that 
\begin{eqnarray*}
\int^\infty_{L}\!\frac{C_1}{|\lambda - i(\omega -R)|} \cdot \frac{1}{1+|\omega |^2}d\omega
 + \int^{-L}_{-\infty}\!\frac{C_1}{|\lambda - i(\omega -R)|} \cdot \frac{1}{1+|\omega |^2}d\omega
< \varepsilon _0.
\end{eqnarray*}
Next, we have
\begin{eqnarray*}
\int^L_{-L}\!\frac{C_1}{|\lambda - i(\omega -R)|} \cdot \frac{1}{1+|\omega |^2}d\omega
 = \int^{R+L}_{R-L}\!\frac{C_1}{|\lambda - i\omega|} \cdot \frac{1}{1+|\omega +R|^2}d\omega,
\end{eqnarray*}
where $1/|\lambda -i\omega |$ is bounded for $-\varepsilon < \lambda < -\varepsilon _0$.
Therefore, for any $\varepsilon _0 > 0$, the estimate
\begin{eqnarray*}
\int^{R+L}_{R-L}\!\frac{C_1}{|\lambda - i\omega|} \cdot \frac{1}{1+|\omega + R|^2}d\omega < \varepsilon _0
\end{eqnarray*}
holds if $R$ is sufficiently large.
This gives
\begin{eqnarray*}
|((\lambda+iR -i\omega )^{-1}\phi, \psi^*)| < 2 \varepsilon _0.
\end{eqnarray*}
It is easy to show that there exists $M_2>0$ such that for any $\varepsilon _0 > 0$, 
\begin{eqnarray*}
|\phi (-i\lambda+R )\psi (-i\lambda+R )g(-i\lambda+R )| < M_2 \varepsilon _0
\end{eqnarray*}
holds if $R$ is sufficiently large.
The other terms included in $\mathcal{R}(\lambda +iR)$ are estimated in the same way.
As a result, it turns out that there exists $M_3 > 0$ such that for any $\varepsilon _0 > 0$,
the inequality $|e^{\lambda t}\mathcal{R}(\lambda +iR)| < M_3 \varepsilon _0$ holds uniformly in 
$-\varepsilon < \lambda < -\varepsilon _0$ if $R$ is sufficiently large.
This result is combined with (\ref{A1}) to yield 
\begin{eqnarray*}
\left| \int^{iR-\varepsilon }_{iR} e^{\lambda t} \mathcal{R}(\lambda ) d\lambda \right|
< (M_3(\varepsilon -\varepsilon _0) + M_1) \varepsilon _0.
\end{eqnarray*}
This proves $\int^{iR-\varepsilon }_{iR} e^{\lambda t} \mathcal{R}(\lambda ) d\lambda \to 0$
as $R\to \infty$.
In a similar manner, we can verify that the other three integrals in (\ref{A0}) tend to zero as $R\to \infty$,
which proves the desired result. $\Box$


\end{document}